# EULER'S GRAPH WORLD
## - MORE CONJECTURES ON GRACEFULNESS BOUNDARIES-II

Suryaprakash Nagoji Rao*
Surusha5152@hotmail.com
May-Jun 2014

**ABSTRACT**. The subclass of Euler graphs with only one type of cycles under (mod 4) operation was studied in Part-1 of this series. It was established that such graphs under regularity are nonexistent for degree >2. Here we consider the subclass of Euler graphs with only two types of cycles under (mod 4) operation. Six cases arise. The case when the cycle types are 0&2(mod 4), the well known class of bipartite Euler graphs, the existence of regular bipartite Euler graphs is very much known. In the other five cases, it transpires that regular Euler graphs with only two types of cycles are nonexistent. Some constructions of Euler graphs with the property are given. We investigate some properties of cycle decompositions, block structure and cycle intersections. Cycle decomposition of Euler graphs allows segregating Euler graphs satisfying Rosa-Golumb criterion and so are nongraceful. In the other case the graphs are possible candidates for gracefulness and are conjectured graceful leading to better understanding of gracefulness boundaries. The cases when the cycle types are 1&2, 1&3, 2&3(mod 4) the graphs are proved to be planar and in other three cases the graphs may not be planar with examples of nonplanar graphs. Probe into the properties of Euler graphs is to propose expected gracefulness boundaries which may guide where to look for graceful graphs. Further, the conjectures may lead to analytical techniques for establishing gracefulness property.

## 1. INTRODUCTION

The word *graph* will mean a finite, undirected graph without loops and multiple edges. Unless otherwise stated a graph is connected. For terminology and notation not defined here we refer to Harary (1972), Mayeda (1972), Buckly and Harary (1988). A graph G is called a *labeled graph* when each node u is assigned a label $\varphi(u)$ and each edge uv is assigned the label $\varphi(uv)=|\varphi(u)-\varphi(v)|$. In this case $\varphi$ is called a *labeling* of G. Define $N(\varphi)=\{n\in\{0,1,...,q_0\}: \varphi(u)=n, \text{ for some } u\in V\}$, $E(\varphi)=\{e\in\{1,2,...,q_0\}: |\varphi(u)-\varphi(v)|=e, \text{ for some edge } uv\in E(G)\}$. Elements of $N(\varphi)$, $(E(\varphi))$ are called *node (edge) labels* of G with respect to $\varphi$. A (p,q)-graph G is *gracefully labeled* if there is a labeling $\varphi$ of G such that $N(\varphi)\subseteq\{0,1,...,q\}$ and $E(\varphi)=\{1,2,...,q\}$. Such a labeling is called a *graceful labeling* of G. A *graceful graph* can be gracefully labeled, otherwise it is a *nongraceful graph* (see Rosa (1967,1991), Golomb (1972), Sheppard (1976) and Guy (1977) for chronology 1969-1977 and Bermond (1978). For a recent survey on graph labelings and varied applications of labeled graphs refer Gallion (2013).

A *(p,q)-graph* has p nodes and q edges; p and q are called *order* and *size* of the graph respectively. A *cycle graph* consists of a single cycle. $C_n$ or *n-cycle* denotes a cycle graph of length n≥3. *Euler trail* is a trail in a graph which visits every edge exactly once. Similarly, Euler circuit is Euler trail which starts and ends on the same node. *Euler graphs* admit Euler circuits and are characterized by the simple criterion that all nodes are of even degree. That is, evenness of node degrees characterizes Euler graphs. *Pendant free graph* G has node degrees two or more. Euler graphs are pendant free graphs. Regular graphs of degree >1 are pendant free graphs. *Core graph* is obtained from G by deleting all pendant nodes recursively till no pendant nodes exist. Core graph of a graph is pendant free. Core graph of a tree is empty graph. Core graph of a unicyclic graph is cycle graph. Core graph of pendant free graph is the graph itself.

By *planting a graph G onto a graph H* we mean identifying a node of G with a node of H. For a pendant free graph G of order p, a *graphforest* GF(G) is constructed from G by planting any number of trees at each node of G. GF(G) is trivial when it is of order p, ie., GF(G)=G and the trees planted are $K_1$s. We assume GF(G) is nontrivial and is of order larger than p. The smallest nontrivial graphforest has pendant node planted at a node. The choice of planting a tree at a node is random. That is, trees of any order and in any number may be planted at each node. The requirement of G a pendant free graph is for convenience. If G is not pendant free then start with core graph of G. Note that for any graph G the simplest graph structure containing G of higher order is graphforest GF(G). A graphforest is *graphtree* when one nontrivial tree is planted. Graphforest of Euler graph G is called *Eulerforest* EF(G). *Cycleforest* CF($C_n$) and *Treeforest* TF(T) are similarly defined. The class of cycleforests is precisely the class of unicyclic graphs. Note that graphforest of a forest is forest and graphforest of a tree is a tree.

---

* Name changed from Suryaprakash Rao Hebbare to Suryaprakash Nagoji Rao, ONGC office order No. MRBC/Expln/1(339)/86-II, Dated 08 March 1999 as per Gazette notification of Govt. Maharashtra December 1998.

**Cycle Decompositions of Euler Graphs**

Euler graph has the novel property that its edge set may be decomposed into edge disjoint cycles and is called *cycle decomposition*. Given a set of cycles there exists Euler graph with this set of cycles as cycle decomposition. Note that, Euler graph, corresponding to a given set of cycles, is not necessarily unique by order. Existence of Euler graph with this set as cycle decomposition follows easily as the graph constructed with these cycles having a common node is Euler. The problem of finding realizable orders or the minimum order for which Euler graphs exist with this set of cycles as cycle decomposition is nontrivial. Cycle decomposition of Euler graph may not be unique.

Some notations follow: For a given cycle decomposition of G classify the cycles $C_n$ into four classes $\mathcal{C}_i$ - the class consisting of cycles $C_n$, $n \equiv i \pmod 4$, $i=0,1,2,3$. $\xi_i = |\mathcal{C}_i|$ - number of cycles in $\mathcal{C}_i$; $C_{ij}$ - $j^{th}$ cycle in the $i^{th}$ class $\mathcal{C}_i$; $s_{ij}$ - length of the cycle $C_{ij}$; $q_i$ - sum of the number of edges in $\mathcal{C}_i$, $i=0,1,2$ or $3$. See Fig.1a. Then for $i^{th}$ class $\mathcal{C}_i$

$$q_i = \sum_{j=1}^{\xi_i} s_{ij} = s_{i1}+s_{i2}+\ldots+s_{i\xi_i}, \quad i=0,1,2,3.$$

By definition, $s_{ij}$ is of the form $4t_{ij}+i$, $i=0,1,2$ or $3$, $t_{ij} \geq 0$. The above identity yields

$$q_i = \sum_{j=1}^{\xi_i}(4t_{ij}+i) = \sum_{j=1}^{\xi_i} 4t_{ij} + \sum_{j=1}^{\xi_i} i,$$

$$= 4\sum_{j=1}^{\xi_i} t_{ij} + i\xi_i,$$

$$\equiv i\xi_i \pmod 4.$$

It follows that, $q = \sum_{i=0}^{3} q_i$ or

$q = q_0+q_1+q_2+q_3$
$\equiv 0\xi_0+1\xi_1+2\xi_2+3\xi_3 \pmod 4$
$\equiv \xi_1+2\xi_2+3\xi_3 \pmod 4 \qquad \ldots(1)$

We summarize the above in

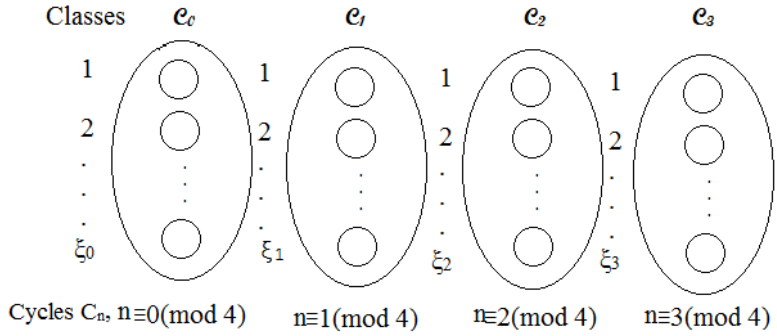

Fig.1a  Typical Cycle Decomposition of Euler Graph.

**Theorem 1**. For any cycle decomposition of Euler (p,q)-graph $q \equiv \xi_1+2\xi_2+3\xi_3 \pmod 4$ holds, where $\xi_i$ is the number of cycles in $\mathcal{C}_i$, $i=0,1,2,3$.

**Corollary 1.1.** For Euler bipartite (p,q)-graph, $q \equiv 2\xi_2 \pmod 4$ holds.

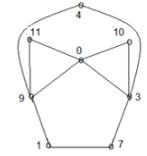
Fig.1b Example of Euler graph with 3,4,5,6,7-cycles.

**EULER GRAPHS**

A number of properties of Euler graphs are known including characterizations (e.g., Zsolt (2010)). Cycle decomposition of Euler graph consists of cycles each edge occurring once. Euler graphs have no pendant nodes and their cycle decomposition may consist of cycles of any length. Fig.1b is a graph with cycles $C_n$, n=3,4,5,6,7. We shall focus on Euler graphs with restricted cycle structure or presence of certain types of cycles and so cycles in its cycle decompositions. Euler graphs exhibit yet more interesting properties so far unknown. As for example, regular Euler graphs of degree >2 with any cycle decomposition having only one class $\mathcal{C}_i$, $i=0,1,2$ or $3$ are nonexistent. These properties and others help in understanding boundaries of gracefulness leading to propose more conjectures.

Some characterizations and properties of Euler graphs (See Zsolt (2010)):

**Theorem A**. The number of edge disjoint paths between any two nodes of Euler graph is even.

**Theorem B** (Euler (1736), Hierholzer (1893)). A connected graph G is Euler graph if and only if all nodes of G are of even degree.

**Theorem C** (Veblen (1912-13)). A connected graph G is Euler if and only if its edge set can be decomposed into cycles.

**Theorem D** (Toida (1973), McKee (1984)). Fleischner (1989,1990).) A connected graph is Euler if and only if each of its edges lies on an odd number of cycles.

**Theorem E** (Bondy and Halberstam (1986)). A graph is Euler if and only if it has an odd number of cycle decompositions.

**Observation 1.** Every cycle decomposition of a bipartite Euler graph has only even cycles.

**Observation 2.** A non-bipartite Euler graph has at least one cycle decomposition consisting of an odd cycle. However, it may have cycle decompositions consisting of only even cycles. For example, Euler graph with a 4-cycle on each edge of a triangle is non-bipartite Euler graph and decomposes into three 4-cycles; or into a 3-cycle and a 9-cycle.

Gracefulness exhibits affinity with Euler graphs as in:

**Theorem F** (Rosa (1967), Golomb (1972)). A necessary condition for Euler $(p,q)$-graph to be graceful is that $[(q+1)/2]$ is even.

Euler graph satisfying this condition is called here a *Rosa-Golomb graph*. This implies that Euler graphs with $q \equiv 1 \text{ or } 2 \pmod 4$ are nongraceful. Most of the known nongraceful graphs contain a subgraph isomorphic to Rosa-Golomb graph. The necessary condition for Euler graphs leads to understanding the intrinsic cycle structure properties of graceful Euler graphs.

## GRACEFULNESS BOUNDARIES THROUGH EULER GRAPHS

### Euler Graphs with Only Two Types of Cycles

A natural generalization of the class of trees viz., graphs with only one type of cycles under (mod 4) operation was studied in Rao (2014). Here, the work is continued for the class of graphs with only two types of cycles in any cycle decomposition of Euler graph. Denote the class of Euler graphs G by $\varepsilon_{ij}$ with only two classes $C_i$ and $C_j$ in every cycle decomposition of G. That is, graphs having only two types of cycles $C_i$ and $C_j$, $i \neq j$, $i,j \equiv 1,2,3 \text{ or } 4 \pmod 4$. This assumption simplifies considerably the general structure of Euler graphs and stand next to trees and graphs with only one type of cycles including unicyclic graphs. Fig.2 shows, Euler and Non-Euler, graphs with types of cycles present. Here, i,j,k,l means cycles in the graph of the type $\equiv i,j,k,l \pmod 4$.

Six cases arise with only two types of cycles under (mod 4) operation: 1) $n \equiv 0 \& 1 \pmod 4$, 2) $n \equiv 0 \& 2 \pmod 4$, 3) $n \equiv 0 \& 3 \pmod 4$, 4) $n \equiv 1 \& 2 \pmod 4$, 5) $n \equiv 1 \& 3 \pmod 4$, 6) $n \equiv 2 \& 3 \pmod 4$. Case-2 corresponds to the class of bipartite Euler graphs. It emerges that regularity of degree greater than two under this condition is impossible within the class of Euler graphs except for the Case 2. We prove that the only regular Euler graphs with just two types of cycles are regular bipartite Euler graphs of degree >2. This characterizes regular bipartite Euler graphs in the class of Euler graphs with only two types of cycles.

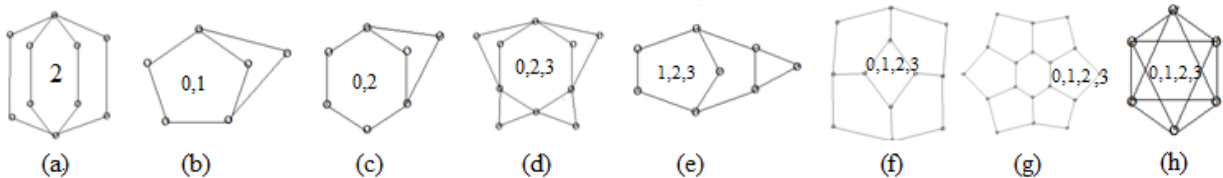

**Fig.2** Graphs not necessarily Euler and Cycle Types present.

**Construction.** $C_n$ be a cycle graph with $n = n_1 n_2$, $n_1 > 1$ and $n_2 > 1$. There are $n_1$ consecutive node pairs at distance $n_2$ with edge disjoint paths. Consider a pair of nodes along $C_n$ at a distance of $n_2$. Add $m > 0$ paths of length $n_2$, $P_{n_2}$ between these $n_1$ pairs of nodes. This results in a $(n + m n_1 (n_2 - 1), n(m+1))$-graph. Node degree is either two or $2(m+1)$ and so the graph is Euler. Each new node added within the paths are of degree 2 and degree of each of $n_1$ nodes is 2(m+1). For q=n(m+1)≡0or3(mod 4) these graphs are candidates for gracefulness. For example, the graph in Fig.2a may be constructed from $C_6$, 6=2x3, by adding two paths $P_3$ at the two nodes at distance 3.

**Cycle Intersection**

Two cycles in a graph may be disjoint or may intersect. The intersection may be a node or a path of length >0. If two cycles have a node in common then no new combined cycle is formed. Hereafter by intersection of two cycles we mean that the cycles intersect in a path of length >0.

For a graph G in $\varepsilon_{ij}$ two intersecting cycles of type i and j with a path $P_t$ of length t>0 and i,j,t=0,1,2,3, the combined cycle will be of type i+j-2t(mod 4) as shown in Table-1. We say that the combined cycle in G is *closed* if it is of type i or j otherwise it is *open*. For example, for i=0 and j=1 the combined cycle for t even is type 1 so closed whereas for t odd is type 3 so open. Table-1 shows closed or open property for values i and j as a row. For a graph G in $\varepsilon_{ij}$ closed or open property is decided for values i or j based on the three rows (i,i), (i,j), (j,j). for example, for a graph in $\varepsilon_{01}$ the combined cycle in the case (1,1) is open but closed as a cycle type 0 results for t odd which is admissible.

| Table-1 Intersection (t) & Combined Cycle | | | | | | | |
|---|---|---|---|---|---|---|---|
| | | | ≡i+j-2t(mod 4) | | | | |
| i | j | i+j | t:0 | 1 | 2 | 3 | Closed for t |
| 0 | 0 | 0 | 0 | 2 | 0 | 2 | Even |
| 0 | 1 | 1 | 1 | 3 | 1 | 3 | Even |
| 0 | 2 | 2 | 2 | 0 | 2 | 0 | Even & Odd |
| 0 | 3 | 3 | 3 | 1 | 3 | 1 | Even |
| 1 | 1 | 2 | 2 | 0 | 2 | 0 | None |
| 1 | 2 | 3 | 3 | 1 | 3 | 1 | Odd |
| 1 | 3 | 4 | 0 | 2 | 0 | 2 | None |
| 2 | 2 | 4 | 0 | 2 | 0 | 2 | Odd |
| 2 | 3 | 5 | 1 | 3 | 1 | 3 | Odd |
| 3 | 3 | 6 | 2 | 0 | 2 | 0 | None |

| Open | | Closed |
|---|---|---|

**Theorem 2.** Let G be a graph in $\varepsilon_{ij}$, i,j=0,1,2,3 which is a block then there are two intersecting cycles of type i and j.

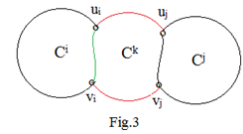

Fig.3

Proof. Since G is from $\varepsilon_{ij}$ it has at least one cycle of type i and type j. Consider $C^i$ and $C^j$ of type i and j respectively, see Fig.3. Since G is a block, for $u_i$ on $C^i$ and $u_j$ on $C^j$ there is a cycle $C^k$ containing $u_i$ and $u_j$. Suppose $C^k$ has $v_i$ and $v_j$ from $C^i$ and $C^j$ in common, respectively. Note that the nodes $u_i, u_j, v_i$ and $v_j$ are distinct otherwise one of them is a cutnode. By definition either $C^k$ is type i or type j. In the first case $C^k$ and $C^j$; and in the second case $C^i$ and $C^k$ are the required intersecting cycles, completing the proof. □

In the ensuing paragraphs we shall consider the six cases separately case by case. In each of these cases three subcases arise for two intersecting cycles and when there is no confusion are denoted by type (i,i) or (j,j) or mixed type (i,j). Further, we denote a block of G with only cycle types ≡i(mod 4), ≡i&j(mod 4) by $B_i$, $B_{ij}$ respectively. Let $\beta_i$, $\beta_{ij}$ denote the number of blocks in G with only cycle types ≡i(mod 4), ≡i&j(mod 4) respectively. Note that, either $\beta_{ij}=0$, $\beta_i>0$ and $\beta_j>0$ holds or $\beta_{ij}>0$ and $\beta_i\geq 0$ or $\beta_j\geq 0$ holds.

**Case-0&1. $\varepsilon_{01}$: Euler graphs with only two types of cycles $C_n$, n≡0&1(mod 4)**

**Construction**. Graphs of the type $\varepsilon_{01}$ may be constructed from G in $\varepsilon_{01}$ by adding any number of paths of same length between two nodes u,v at even distance making sure that every cycle in G containing such a path is of length ≡0or1(mod 4). As an example, for two nodes at even distance of the cycle graph $C_n$, n≡0or1(mod 4) add any number of paths of the same length. The resulting graph belongs to $\varepsilon_{01}$.

**Examples**: Some examples of graphs from $\varepsilon_{01}$ are given in Fig.4. Euler (7,9)-graph (see Fig.4a) with a cycle graph $C_5$ and two nodes adjacent to a pair of nonadjacent nodes of $C_5$ is the smallest graph in $\varepsilon_{01}$. Fig.4a,c are nongraceful. Graceful (15,19)-graph and graceful (16,19)-graph are shown in Fig.4b,d. Note that each block of Fig.4d is a cycle graph $C_n$, $n \equiv 0$ or $1 \pmod 4$.

**Observation 3.** The graph in Fig. 4c, with three 5-cycles having a common edge, belongs to $\varepsilon_{01}$ as it has only 5- and 8-cycles and so cycle types (0,1). The same holds with three intersecting cycles of the type $\equiv 3 \pmod 4$ with a common edge, belongs to $\varepsilon_{03}$ as it has only 7- and 12-cycles and so of the type (0,3). These graphs may be generalised for the classes $\mathcal{G}_1$ and $\mathcal{G}_3$.

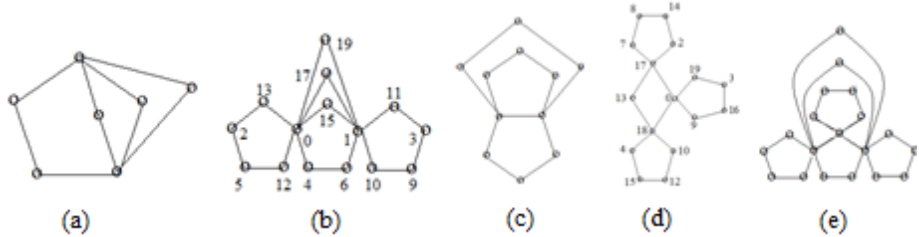

**Fig.4** Graphs from $\varepsilon_{01}$.

| Table-3 Intersection (t) & Combined Cycle | | | | | | | |
|---|---|---|---|---|---|---|---|
| (0,1) Case | | | $\equiv i+j-2t \pmod 4$ | | | | |
| i | j | i+j | t:0 | 1 | 2 | 3 | Closed for t |
| 0 | 0 | 0 | 0 | 2 | 0 | 2 | Even |
| 0 | 1 | 1 | 1 | 3 | 1 | 3 | Even |
| 1 | 1 | 2 | 2 | 0 | 2 | 0 | None |

| Open | | Closed |
|---|---|---|

Table-3 shows combined cycle types for a graph from $\varepsilon_{01}$ for the three cases (0,0), (0,1), (1,1). Closed or open property of combined cycle is decided based on all three cases. Though the case (1,1) results in a combined cycle type 0 for t=1or3(mod 4) so is open but is closed for a graph in $\varepsilon_{01}$.

**Theorem 3.** If two intersecting cycles of a graph in $\varepsilon_{01}$ are of type (0,0) or mixed type (0,1) then the common path $P_t$ has even number of edges. If it is of type (1,1) then $P_t$ has odd number of edges.

Proof. The cases (0,0), (0,1) and (1,1) lead to a combined cycle of length $\equiv -2t \pmod 4$, $\equiv 1-2t \pmod 4$ and $\equiv 2-2t \pmod 4$ respectively. For graphs in $\varepsilon_{01}$, from Table-1, the combined cycle is closed when $t \equiv 0$ or $2 \pmod 4$ in the cases (0,0) and (0,1), and when $t \equiv 1$ or $3 \pmod 4$ in the case (1,1). The result follows. □

**Theorem 4.** Every graph in $\varepsilon_{01}$ has a node of degree 2.

Proof. By Theorem 2 there are two intersecting cycles $C^1$ and $C^2$ of type 0 and type 1 as in Fig.5 with the u-u' path in common. Further, we assume that u,u' nodes are at minimum distance with this property. That is, there is no x-x' subpath of u-u' path such that there is a cycle of type 0 or type 1 intersecting with only x-x' path in common. Let v be a node such that uv is an edge on this u-u' path. Such a node v is assured because the u-u' path of intersection is of even length >0. Claim that v is of degree 2 if not, there is a v-v' path for some node v' different from u and u'. Three possibilities arise, viz., v' is on u-u''(=v')-u' path (Fig. 5a) or u-v'-u' path (Fig. 5b) or v-v'-u' path (Fig. 5i). The first case is not possible as the cycle $C^3$ can neither be of type 0 nor type 1. In the second case, $C^3$ cannot be type 0 and so both $C^3$ and $C^4$ are of type 1. Two subcases arise according as v-u' path is an edge or a path of length 3 or more. If vu' is an edge then consider x as in Fig. 5c. Such a choice of node x follows as v is of even degree. x has x-v' path and by symmetry let v' be on u-v'-u' path as shown. Consider the cycles $C^4$:u,v,x-v'-u and $C^5$: v,u'-v'-x,v. $C^3$ is of type 1 implies that both $C^4$ and $C^5$ cannot be of type 1. Assume $C^4$ is of type 0 then $C^5$ has to be of type 1. Then the v,x-v' path is even. On the other hand, the cycles u-w,v,x-v'-u and v,w-u'-v'-x,v are of type 1 and type 0 respectively and so the path w,v,x-v' is even implying that v,x-v' path is both even and odd, a contradiction.

In the case v-u' path is of length 3 or more let x be a node adjacent to v. Such a node is assured as v is of even degree. Let x be a node adjacent to v. Then five cases arise as shown in Figs. 5d to h. The case d is ruled out as the cycle u,v,x-x'-u''-u cannot be type 0. If it is type 1 then the cycle v,w-u'-x'-x,v has to be only type 0. But then v,x-x' path is even. That it is also odd follows from the fact that v'-v path is odd and v'-v,x-x' is even, a contradiction. The case e is impossible as there will be intersecting cycles of type 0 or 1 with smaller path on u-u' path. The case f is impossible because if the cycle v,x-x'-u'-w,v is type 0 then it is a cycle of type 0 with smaller u-u' path. If it is of type 1 then v,w-u' path is both even and odd, a contradiction. In the case g the cycle v,x-x'-v cannot be type 1 as v,w-u' is shorter path. If it is of type 0 then v,x-x' path is even and odd, a contradiction. The case h is also impossible as the cycle u-x'-x,v,u cannot be type 0. If it is type 1 then v,x-x' path becomes even and odd, a contradiction. Lastly the case in Fig. 5i is impossible because the cycle $C^4$ cannot be type 0 and type 1 otherwise it contradicts that u-u' path is minimum. We conclude that v is of degree 2. □

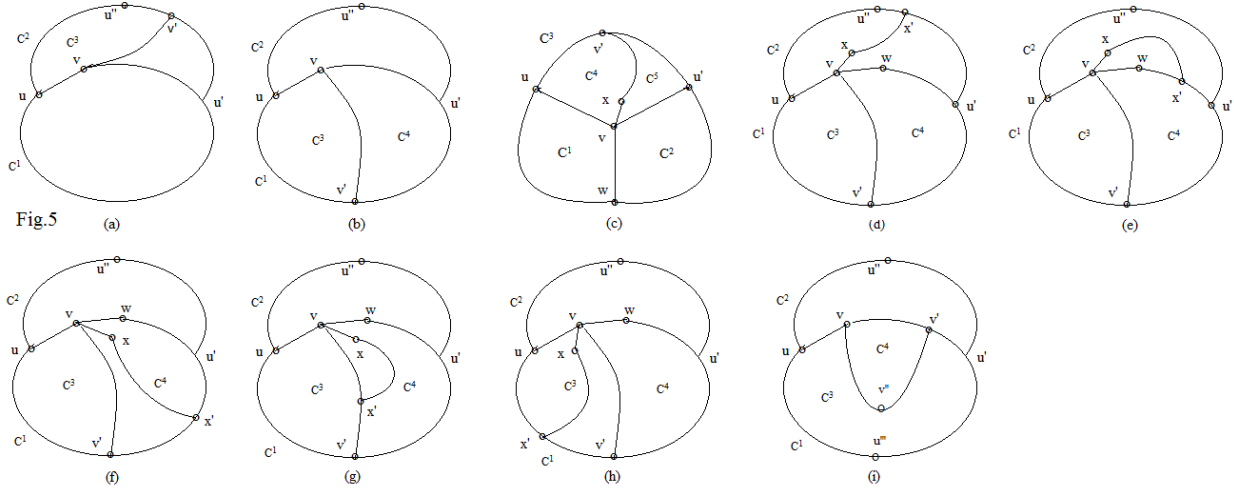

Fig.5 (a) (b) (c) (d) (e) (f) (g) (h) (i)

**Corollary 4.1.** For every graph G in $\varepsilon_{01}$ there is a pair of nodes u,v having exactly two u-v EDPs.

**Corollary 4.2.** Regular Euler graphs in $\varepsilon_{01}$ are nonexistent.

Euler graphs in $\varepsilon_{01}$ may not be planar. Graphs below in Fig. are examples for nonplanar Euler graphs in $\varepsilon_{01}$.

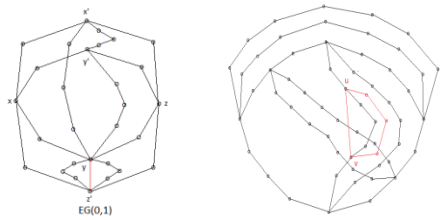

Fig. EG(0,1)

Graphs in $\varepsilon_{01}$, not of Rosa-Golomb type, are candidates for gracefulness. Graphs with each block either of type $\varepsilon_0$ or $\varepsilon_1$, with at least one block of each type, are in $\varepsilon_{01}$. A block itself may be a graph of the type $\varepsilon_{01}$.

**Theorem 5**. A graph G from $\varepsilon_{01}$ has block structure as follows: Either $\beta_{01}>0$ with $\beta_0\geq 0$, $\beta_1\geq 0$ if any from $\varepsilon_0$ or $\varepsilon_1$; or $\beta_{01}=0$ with $\beta_0>0$ and $\beta_1>0$.

**Theorem 6.** Size of a graph from $\varepsilon_{01}$ satisfies $q\equiv\xi_1\pmod{4}$.

Proof follows from Equation (1). Further, if $\xi_1\equiv 1 \text{ or } 2\pmod{4}$ then the graphs are nongraceful; else $\xi_1\equiv 0 \text{ or } 3\pmod{4}$ and the graphs are candidates for gracefulness. Graphs in $\varepsilon_{01}$ satisfy that $\xi_0>0$ and $\xi_1>0$.

**Conjecture 1**. Graphs in $\varepsilon_{01}$ with $\xi_1\equiv 0 \text{ or } 3\pmod{4}$ are graceful.

**Case-0&2.** $\varepsilon_{02}$: **Euler graphs with only two types of cycles $C_n$, n≡0&2(mod 4)**

**Construction**. Graphs of the type $\varepsilon_{02}$ may be constructed from G in $\varepsilon_{02}$ by adding any number of paths of same length between two nodes u,v at even distance making sure that every cycle in G containing such a path is of length ≡0(mod 4). The two graphs in Fig.6a,b,c are of this type. Fig.6d is the hypercube of degree 4, a regular bipartite Euler graph.

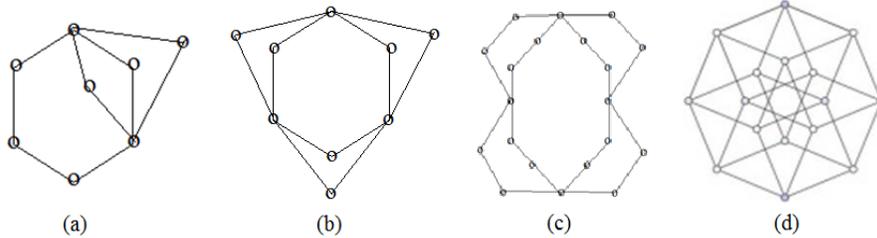

(a)        (b)        (c)        (d)
Fig. 6 Bipartite Graphs from $\varepsilon_{02}$. 4-Cube is shown in Fig. 6(d).

The class of graphs from $\varepsilon_{02}$ is the well known and precisely the class of bipartite Euler graphs. Regularity of any given even degree is realizable. A regular bipartite Euler graph has only even cycles and so are of the type $C_n$, n≡0&2(mod 4). Euler graphs in $\varepsilon_{02}$ may not be planar. Graph in Fig.6d is an example for a nonplanar Euler graph in $\varepsilon_{02}$. Further, Rosa-Golomb criterion for Euler bipartite graphs to be nongraceful reduces to q≡2(mod 4).

| Table-4 Intersection (t) & Combined Cycle | | | | | | | |
|---|---|---|---|---|---|---|---|
| (0,2) Case | | ≡i+j-2t(mod 4) | | | | | |
| i | j | i+j | t:0 | 1 | 2 | 3 | Closed for t |
| 0 | 0 | 0 | 0 | 2 | 0 | 2 | Even & Odd |
| 0 | 2 | 2 | 2 | 0 | 2 | 0 | Even & Odd |
| 2 | 2 | 4 | 0 | 2 | 0 | 2 | Even & Odd |

| Open | | Closed |
|---|---|---|

**Theorem 7.** Irrespective of cycle types, the common path $P_t$ for two intersecting cycles in any graph from $\varepsilon_{02}$ may have even or odd number of edges. In other words, the combined cycle is always even.

**Theorem 8**. A graph G from $\varepsilon_{02}$ has block structure as follows: Either $\beta_{02}>0$ with $\beta_0\geq0$, $\beta_2\geq0$ if any from $\varepsilon_0$ or $\varepsilon_2$; or if $\beta_{02}=0$ then $\beta_0>0$ and $\beta_2>0$.

**Theorem 9.** Size of a graph from $\varepsilon_{02}$ satisfies q≡2$\xi_2$(mod 4), for any $\xi_o>0$, $\xi_2>0$.

Proof follows from Equation (1). By Rosa-Golumb criterion, if $2\xi_2\equiv2$(mod 4) or $\xi_2\equiv1$(mod 2) or $\xi_2$ is odd then the graphs are nongraceful; else $2\xi_2\equiv0$(mod 4) or $\xi_2\equiv0$(mod 2) or $\xi_2$ is even and the graphs are candidates for gracefulness. Graphs in $\varepsilon_{02}$ satisfy that $\xi_0>0$ and $\xi_2>0$.

**Corollary 9.1.** A necessary condition for a graph from $\varepsilon_{02}$ to be graceful is $\xi_2>0$ is even.

**Conjecture 2**. Graphs in $\varepsilon_{02}$ with $\xi_2$ even are graceful.

**Case-0&3.** $\varepsilon_{03}$: **Euler graphs with only two types of cycles $C_n$, n≡0&3(mod 4)**

**Construction**. Graphs of the type $\varepsilon_{03}$ may be constructed from a cycle graph $C_n$, n≡3(mod 4) by adding any number of paths of same length between two nodes u,v at even distance making sure that every cycle in G containing such a path is of length ≡0(mod 4).

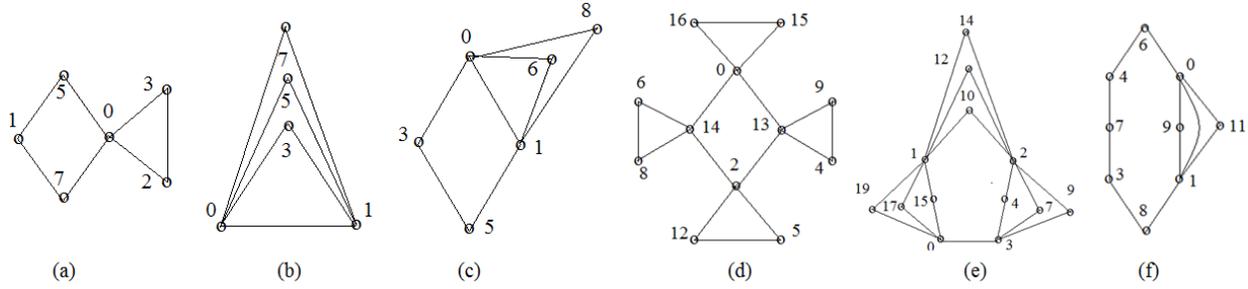

**Fig.7** Graceful Graphs from $\varepsilon_{03}$.

| Table-5 Intersection (t) & Combined Cycle ||||||||
| (0,3) Case ||| $\equiv$i+j-2t(mod 4) ||||
| i | j | i+j | t:0 | 1 | 2 | 3 | Closed for t |
| 0 | 0 | 0 | 0 | 2 | 0 | 2 | Even |
| 0 | 3 | 3 | 3 | 1 | 3 | 1 | Even |
| 3 | 3 | 6 | 2 | 0 | 2 | 0 | Odd |

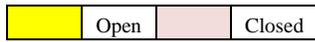

**Theorem 10.** If two intersecting cycles of a graph in $\varepsilon_{03}$ are of type (0,0) or mixed type (0,3) then the common path $P_t$ has even number of edges. If it is of type (3,3) then $P_t$ has odd number of edges.

Proof follows on the similar lines that of Theorem 3.

**Theorem 11.** Every graph in $\varepsilon_{03}$ has a node of degree 2.

The proof of this theorem runs on the similar lines of Theorem 4 in the case $\varepsilon_{01}$ because the rules for the intersections of cycle types (0,0), (0,3) and (3,3) for a graph in $\varepsilon_{03}$ are same as that of (0,0), (0,1) and (1,1) for a graph in $\varepsilon_{01}$.

**Corollary 11.1.** Regular Euler graphs in $\varepsilon_{03}$ are nonexistent.

**Corollary 11.2.** For any graph in $\varepsilon_{03}$ there exists a pair of nodes u,v with exactly two u-v EDPs.

Euler graphs in $\varepsilon_{03}$ may not be planar. Graph in Fig.7a is an example for a nonplanar Euler graph in $\varepsilon_{03}$.

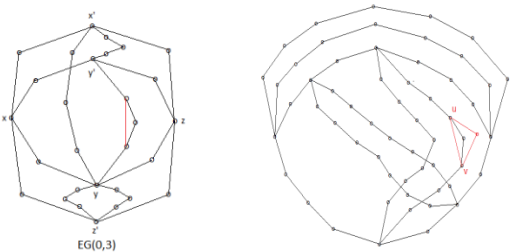

Fig.7a EG(0,3)

Complete graph of order 4 is a regular graph of degree 3 and so not Euler. Further, it has cycle types 0 or 3 only.

**Theorem 12.** A graph G from $\varepsilon_{03}$ has block structure as follows: Either $\beta_{03}>0$ with $\beta_0\geq 0$, $\beta_3\geq 0$ if any from $\varepsilon_0$ or $\varepsilon_3$; or if $\beta_{03}=0$ then $\beta_0>0$ and $\beta_3>0$.

**Theorem 13.** Size of a graph from $\varepsilon_{03}$ satisfies, $q\equiv 3\xi_3\pmod{4}$, for any $\xi_0>0$ and $\xi_3>0$.

Proof follows from Equation (1). By Rosa-Golumb criterion, if $3\xi_3\equiv 1\,or\,2\pmod{4}$ or $\xi_3\equiv 2\,or\,3\pmod{4}$ then the graphs are nongraceful; else $3\xi_3\equiv 0$ or $3 \pmod 4$ or $\xi_3\equiv 0$ or $1 \pmod 4$ and the graphs are candidates for gracefulness. Graphs in $\varepsilon_{03}$ satisfy that $\xi_0>0$ and $\xi_3>0$.

**Conjecture 3.** Graphs in $\varepsilon_{03}$ with $\xi_3\equiv 0$ or $1 \pmod 4$ are graceful.

**Case 1&2.** $\varepsilon_{12}$: Euler graphs with only two types of cycles $C_n$, $n\equiv 1\&2 \pmod 4$

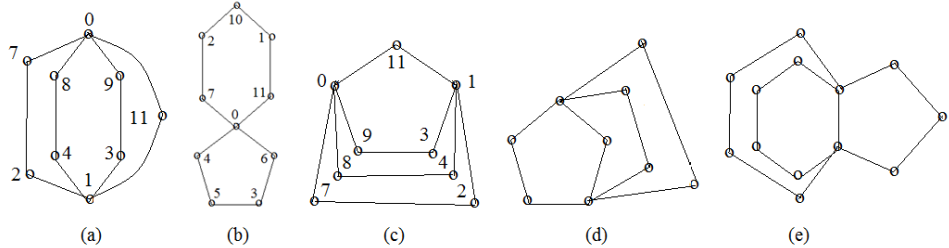

(a) (b) (c) (d) (e)

**Fig.8** Graphs from $\varepsilon_{12}$ with graceful numbering for the first three.

**Theorem 14.** If two intersecting cycles of a graph in $\varepsilon_{12}$ are of type (2,2) or mixed type (1,2) then the common path $P_t$ has odd number of edges. If it is of type (1,1) then $P_t$ has even number of edges.

| Table-6 Intersection (t) & Combined Cycle | | | | | | | |
|---|---|---|---|---|---|---|---|
| (1,2) Case | | | $\equiv i+j-2t \pmod 4$ | | | | |
| i | j | i+j | t:0 | 1 | 2 | 3 | Closed for t |
| 1 | 1 | 2 | 2 | 0 | 2 | 0 | None |
| 1 | 2 | 3 | 3 | 1 | 3 | 1 | Odd |
| 2 | 2 | 4 | 0 | 2 | 0 | 2 | Odd |

| | Open | | Closed |
|---|---|---|---|

**Theorem 15.** If cycles of type 1 ($C^1$) and type 2 intersect then no type 1 cycle neither contains nor contained in type 1 cycle $C^1$ with the same intersection.

The result follows as the combined cycle formed out of any two cycles of type 1 is of type 0.

**Theorem 16.** The graphs in Fig.9 are forbidden in $\varepsilon_{12}$.

**Proof.** The proof technique we follow is to prove the existence of a path which is both even and odd in length. We infer from Table-6 that the path of intersection is even for cycles of type (1,1) and odd for cycles of type (1,2) and (2,2). Further, the combined cycle in cases (1,1) and (2,2) is type 2 while in case (1,2) is type 1. Now, consider the subgraph in Fig.9a. The w-v and v-x paths are odd because they are intersection of cycles type 2 and type 1. So, w-v-x path is of even length. However, w-v-x is intersection of cycles type 2 and type 1 and so is of odd length. Therefore, length of w-v-x path is both even and odd, a contradiction and so is forbidden as subgraph. On similar lines it may be proved that the subgraphs in Figs.9 b to g are forbidden by considering the paths w-v-x in Figs. 9b,c,d; w-x path in Figs.9e,f; and u-w-x' path in Fig.9g. □

**Theorem 17.** A graph G from $\varepsilon_{12}$ has block structure as follows: Either $\beta_{12}>0$ with $\beta_1\geq 0$, $\beta_2\geq 0$; or if $\beta_{12}=0$ with $\beta_1>0$ and $\beta_2>0$.

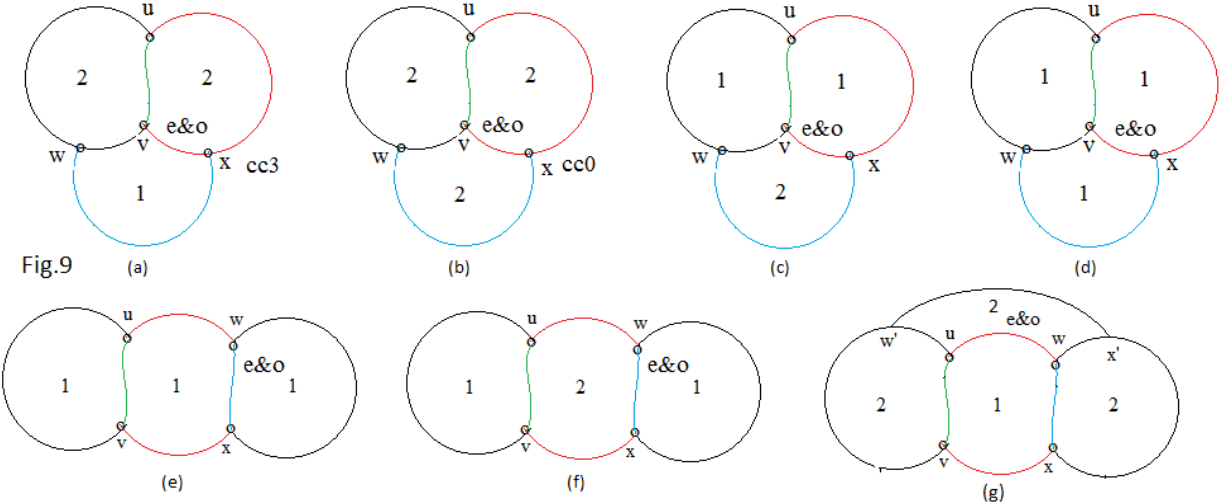

Fig.9 (a) (b) (c) (d) (e) (f) (g)

**Theorem 18.** Every graph in $\varepsilon_{12}$ has a node of degree 2.

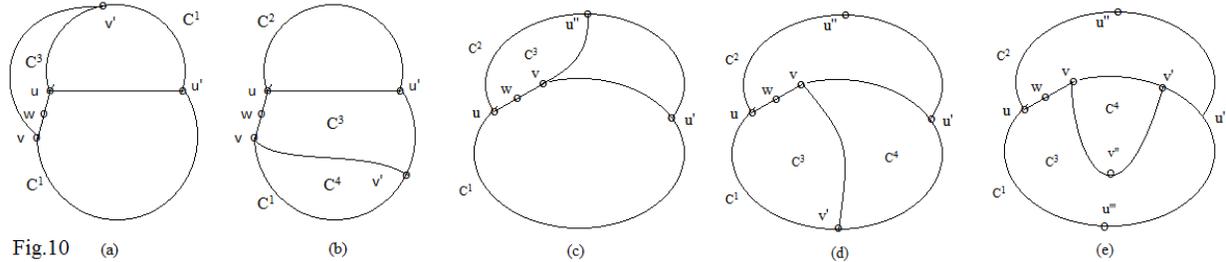

Fig.10 (a) (b) (c) (d) (e)

**Proof.** By Theorem 2 there are two intersecting cycles $C^1$ and $C^2$ of type 1 and type 2 as in Fig.10. Further, we assume that there is no x-x' path on u-u' path such that there is a cycle of type 1 or type 2 intersecting in only x-x' path. In other words, the nodes u,u' are so chosen that they are at minimum distance. Consider the u-u' path of intersection of the cycles $C^1$ and $C^2$ and is of odd length. Two cases arise according as the odd u-u' path is an edge or a path of length 3 or more. In the first case, if uu' is an edge then there exist nodes w,v such that v≠u' and uw, wv are edges on $C^1$. Existence of such nodes v,w follows from cycles of type 2. Further, we assume that $C^1$ is smallest containing v,w,u-u' path. Claim that v is of degree 2. If not then there is a v-v' path and two possibilities arise as in Figs. 10a,b. In the case a, the cycle $C^3$ can neither be type 2 nor type 1. The case b is impossible by the choice of $C^1$.

In the second case, suppose u-u' path is of length 3 or more with nodes v,w such that uw, wv are edges. Claim that v is of degree 2. If not, there exists a third node v' and a new v-v' path. Three possibilities arise, viz., v' is on u-u''-u' path or u-v'-u' path or v-v'-u' path as shown in Figs. 10c,d,e. In the case c, the cycle $C^3$ can neither be of type 2 nor type 1. In the case d, the cycle $C^3$ cannot be of type 2. If $C^3$ is of type 1 then $C^4$ can be only of type 1. But then, the v-u' path is even, a contradiction as u-u' path and so v-u' path is also odd, a contradiction. Lastly, the case e is ruled out otherwise it contradicts that the u-u' path is minimum. We conclude that v is of degree 2 completing the proof. □

**Corollary 18.1.** For any graph in $\varepsilon_{12}$ there exists a pair of nodes u,v with exactly two u-v EDPs.

**Corollary 18.2.** Regular Euler graphs in $\varepsilon_{12}$ are nonexistent.

**Theorem** 19a. Euler graphs in $\varepsilon_{12}$ are planar.

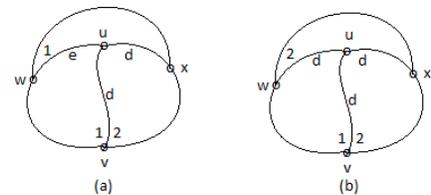

(a) (b)

**Proof.** By Theorem ... there are two intersecting cycles C1 and C2. Fig a shows intersecting cycles =1 and =2 (mod 4). the length of intersection is odd. the combined cycle of C1 and C2 be C3 which is =1(mod 4). let w and x be vertices as shown and w-x be a path apart from the path along C3 forming a cycle C4. Two possibilities arise as C4 is =1 or 2(mod 4) as shown in Fig a,b. In first case, C1, C4 intersect in an even path and C2, C4 intersect in an odd path. however, C3 and C4 intersect in even path which is contrary to the path of intersection is odd. In the second case, C1 and C2 intersect in odd paths with C4 implying w-u-x path is even. on the other hand, C3 and C4 intersect in odd path leading to contradiction.

**Theorem 19**. Size of a graph from $\varepsilon_{12}$ satisfies, $q \equiv \xi_1 + 2\xi_2 \pmod 4$, for any $\xi_1 > 0$ and $\xi_2 > 0$.

Proof follows from Equation (1). By Rosa-Golumb criterion, if $\xi_1 + 2\xi_2 \equiv 1 \text{ or } 2 \pmod 4$ then the graphs are nongraceful; else $\xi_1 + 2\xi_2 \equiv 0 \text{ or } 3 \pmod 4$ and the graphs are possible candidates for gracefulness. Graphs in $\varepsilon_{12}$ satisfy that $\xi_1 > 0$ and $\xi_2 > 0$.

**Corollary 19.1.** Two necessary conditions for a graph from $\varepsilon_{12}$ is graceful are: If $\xi_1 + 2\xi_2 \equiv 0 \pmod 4$ then $\xi_1$ is even; else $\xi_1 + 2\xi_2 \equiv 3 \pmod 4$ then $3 - \xi_1$ is even implying that $\xi_1$ is odd.

**Conjecture 4**. Graphs in $\varepsilon_{12}$ with $\xi_1 + 2\xi_2 \equiv 0 \text{ or } 3 \pmod 4$ are graceful.

**Case 1&3. $\varepsilon_{13}$: Euler graphs with only two types of cycles $C_n$, $n \equiv 1 \& 3 \pmod 4$**

This class of graphs has simple structure.

| Table-7 Intersection (t) & Combined Cycle | | | | | | |
|---|---|---|---|---|---|---|
| (1,3) Case | | $\equiv i+j-2t \pmod 4$ | | | | |
| I | j | i+j | t:0 | 1 | 2 | 3 | Closed for t |
| 1 | 1 | 2 | 2 | 0 | 2 | 0 | None |
| 1 | 3 | 4 | 0 | 2 | 0 | 2 | None |
| 3 | 3 | 6 | 2 | 0 | 2 | 0 | None |

| | Open | | Closed |
|---|---|---|---|

It clearly follows from Table-7 that the combined cycle is open in all three cases irrespective of t is even or odd.

**Theorem 20.** No two cycles of a graph in $\varepsilon_{13}$ intersect in a path $P_t$, $t > 0$.

**Proof.** For a graph in $\varepsilon_{13}$ if there exist two cycles of type i and j which intersect in a path of length $t > 0$ then the combined cycle is of length $\equiv i+j-2t \pmod 4$. This for $i=1$ and $j=3$ reduces to $\equiv 4-2t \pmod 4$. For t even or odd the combined cycle is $\equiv 0 \text{ or } 2 \pmod 4$, a contradiction. □

**Corollary 20.1.** Every block of a graph in $\varepsilon_{13}$ is a cycle graph.

**Corollary 20.2.** A graph in $\varepsilon_{13}$ has a node of degree 2.

**Proof**. Let G be a graph from $\varepsilon_{13}$. By Theorem 20 if two cycles intersect then they have a node in common and is a cutnode. If G has no node of degree 2 then each node is of degree at least 4. It follows then that G is an infinite graph, a contradiction. □

**Corollary 20.3.** Regular Euler graphs in $\varepsilon_{13}$ are nonexistent.

**Corollary 20.4.** Euler graphs in $\varepsilon_{13}$ are planar.

**Theorem 21**. A graph G from $\varepsilon_{13}$ has block structure as follows: $\beta_{13}=0$ with $\beta_1>0$ and $\beta_3>0$ with each block being a cycle graph $C_n$, $n \equiv 1 \text{ or } 3 \pmod 4$.

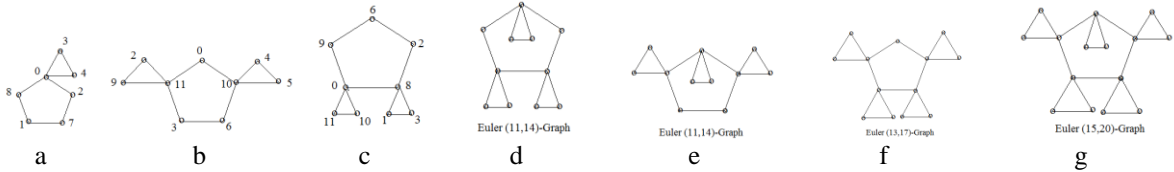

| a | b | c | d | e | f | g |

Fig.11a,b,c Euler Graceful Graphs.   Fig.11d,e: Euler Nongraceful Graphs.   Fig.11f,g: Gracefulness not known.

**Theorem 22.** Size of a graph from $\varepsilon_{13}$ satisfies, $q \equiv \xi_1 + 3\xi_3 \pmod 4$, for any $\xi_1 > 0, \xi_3 > 0$.

Proof follows from Equation (1). By Rosa-Golumb criterion if $\xi_1 + 3\xi_3 \equiv 1 \text{ or } 2 \pmod 4$ then the graphs are nongraceful; else $\xi_1 + 3\xi_3 \equiv 0 \text{ or } 3 \pmod 4$ and the graphs are possible candidates for gracefulness. Graphs in $\varepsilon_{01}$ satisfy that $\xi_1 > 0$ and $\xi_3 > 0$.

**Corollary 22.1.** Two necessary conditions for a graph from $\varepsilon_{13}$ to be graceful are: If $\xi_1 + 3\xi_3 \equiv 0 \pmod 4$ then $3/\xi_1$; else $\xi_1 + 3\xi_3 \equiv 3 \pmod 4$ then $3/(3 - \xi_1)$ which also implies $3/\xi_1$.

**Conjecture 5**. Graphs in $\varepsilon_{13}$ with $\xi_1 + 3\xi_3 \equiv 0 \text{ or } 3 \pmod 4$ are graceful.

**Case 2&3. $\varepsilon_{23}$: Euler graphs with only two types of cycles $C_n, n \equiv 2 \& 3 \pmod 4$**

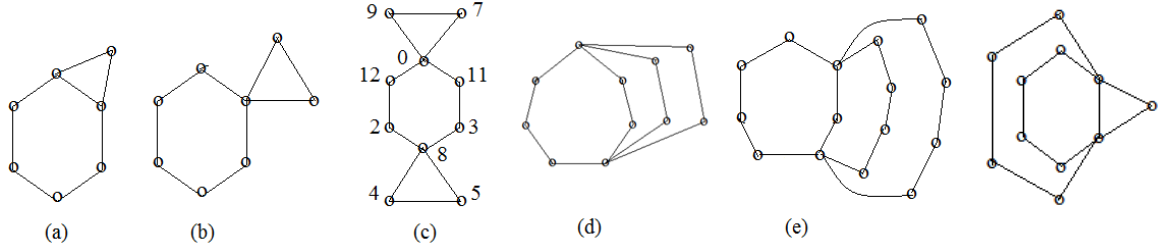

(a)   (b)   (c)   (d)   (e)

Fig.12 Graphs from $\varepsilon_{23}$

**Construction**. Graphs of the type $\varepsilon_{23}$ may be constructed from a cycle graph $C_n, n \equiv 3 \pmod 4$ by adding any number of paths of same length between two nodes u,v at odd distance making sure that every cycle in G containing such a path is of length $\equiv 2 \pmod 4$. An example is given in Fig.12d. Another construction is a graph having at least two blocks with at least one block of the type $\varepsilon_3$ or $\varepsilon_2$. Examples are Fig.12b,c. The graph in Fig.12c is graceful.

| Table-8 Intersection (t) & Combined Cycle | | | | | | |
|---|---|---|---|---|---|---|
| (2,3) Case | | | $\equiv i + j - 2t \pmod 4$ | | | |
| i | j | i+j | t:0 | 1 | 2 | 3 | Closed for t |
| 2 | 2 | 4 | 0 | 2 | 0 | 2 | Odd |
| 2 | 3 | 5 | 1 | 3 | 1 | 3 | Odd |
| 3 | 3 | 6 | 2 | 0 | 2 | 0 | Even |

| | Open | | Closed |
|---|---|---|---|

**Theorem 23.** If two intersecting cycles of a graph in $\varepsilon_{23}$ are of type (2,2) or mixed type (2,3) then the common path $P_t$ has odd number of edges. If it is of the type (3,3) then $P_t$ has even number of edges.

**Theorem 24**. A graph G from $\varepsilon_{23}$ has block structure as follows: Either $\beta_{23} > 0$ with $\beta_2 \geq 0, \beta_3 \geq 0$ if any from $\varepsilon_2$ or $\varepsilon_3$; or if $\beta_{23} = 0$ then $\beta_2 > 0, \beta_3 > 0$.

**Theorem 25.** Every graph in $\varepsilon_{23}$ has a node of degree 2.

The proof of this theorem runs on the similar lines of Theorem 18 in the case $\varepsilon_{12}$ because the rules for the intersections of cycle types (2,2), (2,3) and (3,3) for a graph in $\varepsilon_{32}$ are same as that of (2,2), (2,1) and (1,1) for a graph in $\varepsilon_{12}$.

**Corollary 25.1.** For any graph in $\varepsilon_{23}$ there exists a pair of nodes u,v with exactly two u-v EDPs.

**Corollary 25.2.** Regular Euler graphs in $\varepsilon_{23}$ are nonexistent.

**Theorem 26a.** Euler graphs in $\varepsilon_{23}$ are planar.

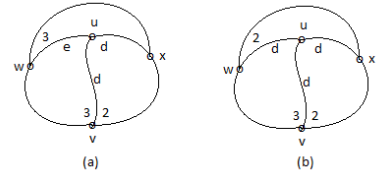

Proof follows similar to the proof of Theorem 19a. See also Figs. a,b.

**Theorem 26**. Size of a graph from $\varepsilon_{23}$ satisfies, $q \equiv 2\xi_2 + 3\xi_3 \pmod{4}$, for any $\xi_2 \geq 0$, $\xi_3 \geq 0$.

**Proof** follows from Equation (1). If $2\xi_2 + 3\xi_3 \equiv 1 \text{ or } 2 \pmod{4}$ then the graphs are nongraceful. Otherwise $2\xi_2 + 3\xi_3 \equiv 0 \text{ or } 3 \pmod{4}$ and the graphs are candidates for gracefulness. Graphs in $\varepsilon_{23}$ satisfy that $\xi_2 > 0$ and $\xi_3 > 0$.

**Corollary 26.1.** Two necessary conditions for a graph from $\varepsilon_{23}$ are: If $2\xi_2 + 3\xi_3 \equiv 0 \pmod{4}$ then $\xi_3$ is even or $3/\xi_2$; else $2\xi_2 + 3\xi_3 \equiv 3 \pmod{4}$ then $1-\xi_3$ is even implying that $\xi_3$ is odd or $3/\xi_2$.

**Conjecture 6.** Graphs in $\varepsilon_{23}$ with $2\xi_2 + 3\xi_3 \equiv 0 \text{ or } 3 \pmod{4}$ are graceful.

**Summary of Six Cases**

It was proved in Rao(2014) Part-I that regular Euler graphs of degree >2, with exactly one class in any cycle decomposition or with only one type of cycles, are nonexistent. As a corollary it follows that a regular Euler graph of degree >2 has two or more cycle types. We now summarize the above six cases:

**Theorem 27.** Euler graph with only two types of cycles is regular iff it is regular bipartite Euler graph of degree >2.

Combining this result with that of Rao(2014) Part-I, Euler graphs with only one or two types of cycles can be characterized:

**Theorem 28.** Euler graph with utmost two types of cycles is regular iff it is an odd cycle graph or a regular bipartite Euler graph.

Lastly, we make the following in support of the conjecture that *planting a tree on a node of a graceful graph continues to remain graceful*, that is, *graphforests of a graceful graph are graceful*.

**Conjecture 7.** Eulerforests of a graph in Conjectures 1 to 6 above are graceful.

**Problems**. Study extremal graphs in $\varepsilon_{ij}$, i,j=0,1,2,3 as a function of size q for given order p. Study the graphs in $\varepsilon_{ij}$, i,j=0,1,2,3 under planarity and corresponding extremal graphs. Establish gracefulness for infinite class of graphs of this type in support of the conjectures.

**Acknowledgements**

The author, an Oil & Gas Professional with works in Graph Theory, expresses deep felt gratitude to
The *Sonangol Pesquisa e Produção*, Luanda, Angola
for excellent facilities, support and encouragement (2011-'13) with special thanks to
Mr. *Joao Noguiera*, Managing Director, Development Directorate.

Deep felt gratitude to *Alexander Rosa*, Emeritus Professor,
Dept. of Mathematics & Statistics, McMaster University, West Hamilton, Ontario, Canada and
*Solomon W Golumb*, an American Mathematician, Engineer and Professor of Electrical Engineering
at the University of Southern California (A recipient of the USC Presidential Medallion,
the IEEE Shannon Award of the Information Theory Society and three honorary doctorate degrees,
National Medal of Science presented by President Barack Obama)
for their initial inspiring works on Graceful Labelings and to
Emeritus Professor *G.A. Patwardhan*, Combinatorics, Department of Mathematics, IIT Bombay.